\documentclass[graybox]{svmult}
\usepackage{amssymb,amsmath}

\usepackage{mathptmx}       
\usepackage{helvet}         
\usepackage{courier}        
\usepackage{type1cm}        
\usepackage{makeidx}         
\usepackage{graphicx}        
\usepackage{multicol}        
\usepackage[bottom]{footmisc}
\usepackage{shuffle}
\input xy
\xyoption{all}

\def\cala{{\cal A}}

\def\E{{\mathcal E}}
\def\K{{\mathbb K}}

\def\Z{{\mathbb Z}}
\def\g{{\mathfrak g}}

\newcommand{\ba}{\begin{array}}
\newcommand{\ea}{\end{array}}
\newcommand{\beq}{\begin{equation}}
\newcommand{\eeq}{\end{equation}}
\newcommand{\beqa}{\begin{eqnarray}}
\newcommand{\eeqa}{\end{eqnarray}}

\begin{document}
\title*{Young tableaux  and   homotopy commutative algebras}

\author{Michel Dubois-Violette and Todor Popov}
\institute{Michel Dubois-Violette \at Laboratoire de Physique Th\'eorique, UMR 8627, Universit\'e Paris XI,\\ B\^atiment 210,
F-91 405 Orsay Cedex, France 
\and Todor Popov \at  Institute for Nuclear Research and Nuclear Energy, Bulgarian Academy of Sciences,\\
72 Tsarigradsko chauss\'ee, 1784 Sofia, Bulgaria,  \email{tpopov@inrne.bas.bg}
\and LPT-ORSAY 12-05, 
Talk given by T.P. at the International Workshop "Lie Theory-9", Varna, 2011 }
\maketitle
\abstract
{A   homotopy commutative algebra, or  $C_{\infty}$-algebra,
 is defined  via the Tornike Kadeishvili homotopy transfer theorem  on the vector space generated by the set of Young tableaux with  self-conjugated Young diagrams $\left\{\lambda:\lambda=\lambda'\right\}$. We prove that this  $C_{\infty}$-algebra 
is generated in degree 1 by the binary and the ternary operations.}

\section{Introduction}
\label{sec:1}

We consider the 2-nilpotent graded Lie algebra $\mathfrak{g}$, with degree one generators in the finite dimensional vector  space $V$ over   a 
field $\K$ of
characteristic 0,
$$
\mathfrak{g} = V \oplus [ V, V ]  \ .
$$
The Universal Enveloping Algebra (UEA) $U\g$ arises naturally in physics as the subalgebra closed by the creation operators of the parastatistics algebra. The algebra of  creation and annihilation parastatistics operators was introduced by 
H.S.\ Green\cite{Green}, its defining relations generalize the canonical (anti)commutation relations.

As an UEA of a finite dimensional positively 
 graded Lie algebra,  $U\g$  belongs to the class of Artin-Schelter regular algebras(see e.g.\cite{FV}).
As every 
finitely generated graded connected algebra,  $U\g$ has a free minimal resolution which is canonically built from the data of its Yoneda algebra 
$ \E := {\rm{ Ext}}_{U\g} (\K, \K).$
By construction the Yoneda algebra $\E$ is isomorphic (as algebra) to the cohomology 
of the Lie algebra (with coefficients in the trivial representation provided by the ground field $\K$)
\beq
\label{iso}
\E =  {{\rm{ Ext}}_{U\g}^{\bullet} (\K, \K)} \cong H^{\bullet} (\mathfrak g  ,  \K) 
\eeq
 the product on $\E$ being the super-commutative wedge product between cohomological classes in $H^{\bullet} (\mathfrak g  ,  \K)$.
 
 An important  
 result due to  J\'ozefiak and Weyman \cite{JW} implies that  a basis of the cohomology  
 $\E=H^{\bullet} (\mathfrak g  ,  \K)$ is indexed by  Young tableaux with
 self-conjugated Young  diagrams (i.e., symmetric with respect to the diagonal). On the other hand according to the homotopy transfer theorem due to  Tornike Kadeishvili \cite{Kadeishvili}  the Yoneda algebra $\E$ is a 
$C_{\infty}$-algebra.
 
 The aim of this note is to endow the cohomology $H^{\bullet}(\mathfrak g  , \K)$
(i.e., the vector space generated by the  set of Young tableaux with  self-conjugated Young diagrams $\left\{\lambda:\lambda=\lambda'\right\}$)
  with a $C_{\infty}$-structure, induced by the
 isomorphism (\ref{iso}) through  the homotopy transfer.
 
Here we deal only with the parafermionic case corresponding to  an (even) vector space $V$.
To include the parabosonic degrees of freedom one have to consider $V$ in the 
 category of vector superspaces. The supercase will be consider elsewhere.


\section{  Artin-Schelter regularity }
\label{sec:2}

Let $\g$ be the 2-nilpotent graded  Lie algebra $\g=V\oplus \bigwedge^2 V $ generated by the finite dimensional  vector  space 
$V$ having   Lie bracket 
   \beq [x,y]:= \left\{ 
   \ba{cccc}
   &x\wedge y &\quad & x,y \in V
   \\& 0 && \quad  \mbox{otherwise} \ea 
   \right.   \ . 
   \eeq
We denote the Universal Enveloping Algebra $U\g$ 
  by $PS$
and will 
refer
to it as {\it parastatistics algebra }(by some abuse\footnote{Strictly speaking $PS(V)$ is the creation 
{ parastatistics algebra}, closed by creation operators alone.}). The parastatistics  algebra $PS(V)$
generated in $V$ is graded 
$$
PS(V):= U\g  = U(V \oplus {\bigwedge}^2 V) = T(V)/ ([[V,V],V]) \ .
$$
We shall write simply $PS$ when the space of generators $V$ is clear from the context.

Artin and Schelter \cite{AS} introduced a class of regular algebras sharing some ``good'' homological properties with the
polynomial algebra $\K[V]$.  These algebras were dubbed Artin-Schelter regular algebras
(AS-regular algebra for short).

\begin{definition}(AS-regular algebras) A connected graded algebra $\cala = \K \oplus \cala_1 \oplus \cala_2 \oplus \ldots$
is called Artin-Schelter regular of dimension $d$ if

(i) $\cala$ has finite global dimension $d$,

(ii) $\cala$ has finite Gelfand-Kirillov dimension,

(iii) $\cala$ is Gorenstein, i.e., ${\rm Ext}^i_{\cala}(\K, \cala)= \delta^{i, d} \K$.
\end{definition}

A general theorem claims that the UEA of a finite dimensional positively graded Lie algebra is an AS-regular algebra 
of global dimension equal to the dimension of the Lie algebra\cite{FV}.
Hence the  parastatistics algebra $PS$ is  AS-regular of global dimension $d=\frac{\dim V(\dim V+1)}{2}$.
In particular the finite global dimension of $PS$ implies that  the ground field $\K$ 
has a  minimal resolution ${P}_{\bullet} $ by projective left $PS$-modules $P_n$
  \beq
  \label{minP}
 { P}_{\bullet}: \qquad  0 \rightarrow P_d \rightarrow  \cdots \rightarrow P_n\rightarrow \cdots \rightarrow P_2 \rightarrow P_1 \rightarrow P_0 \stackrel{\epsilon}{\rightarrow} \K \rightarrow 0  
 \ .
  \eeq
 Here $\K$ is  a trivial left $PS$-module, the action being defined  by the projection $\epsilon$ onto $PS_0=\K$. 
  Since $PS$ is graded and in the category of graded modules projective module
  is the same as free module\cite{Cartan}, we have
  $P_n \cong PS \otimes E_n $ where $E_n $ are finite dimensional vector spaces.
 
 The minimal projective resolution is unique (up to a isomorphism).
 Minimality implies that the  complex $ \K \otimes_{PS} { P}_{\bullet}$ has ``zero differentials''
  hence  
 $$H_{\bullet} (\K \otimes_{PS} { P}_{\bullet}) =\K \otimes_{PS} { P}_{\bullet} =E_n \ .$$ 
 One can calculate the derived functor $\mathrm{Tor}^{PS}_{n}(\K,\K)$   using the resolution $P_{\bullet}$, it yields 
   \beq
   \mathrm{Tor}^{PS}_{n}(\K,\K) =E_n  \ .
   \eeq
 The data of a minimal resolution of $\K$ by free   $PS$-modules provides an easy way to find
 $\mathrm{Tor}^{PS}_{n}(\K,\K)$. Conversely if the  spaces $\mathrm{Tor}^{PS}_{n}(\K,\K)$
 are known then one can construct a  minimal free resolution  of $\K$.

The Gorenstein property 
guarantees
 that when applying the functor ${\rm Hom}_{PS}(-, PS)$ to the minimal free resolution $P_{\bullet}$ we get another minimal free resolution $P^{\bullet}:={\rm Hom}_{PS}({ P}_{\bullet}, PS)$ of $\K$ by right $PS$-modules
 \beq
P^{\bullet}: \quad  0 \leftarrow \K \leftarrow P^{\ '}_d \leftarrow  \cdots \leftarrow P_n^{\ '} \leftarrow \cdots \leftarrow P_2^{\ '} \leftarrow P_1^{\ '} \leftarrow P_0^{\ '} 
 \leftarrow 0  \ 
 \eeq
with $P_n^{\ '}  \cong E_n^{\ast} \otimes PS$. Note that by construction $E_n^{\ast} ={\rm Ext}^n_{PS}(\K,\K)$
thus one has vector space isomorphisms \cite{Cartan}
\beq
\label{cartaniso}
E_n \cong E^\ast_n \cong \mathrm{Tor}^{PS}_{n}(\K,\K) \cong {\rm Ext}^n_{PS}(\K,\K) \ .
\eeq
The Gorenstein property is the analog of the Poincar\'e duality since it implies 
$$
E^{\ast}_{d-n}\cong E_n  \ .
$$
The finite global dimension $d$ of $PS$ and the Gorenstein condition 
imply
 that its Yoneda algebra 
$$
\E^{\bullet}:={\rm Ext}^{\bullet}_{PS}(\K,\K)\cong \bigoplus_{n=0}^d E^{\ast}_n \ 
$$
is Frobenius\cite{LPWZ}.

\section{Homology and cohomology of $\g$ }

 Let us first recall that the standard Chevalley-Eilenberg chain complex $C_{\bullet}(\g)=( U \g  \otimes_{\K} \wedge^p \g , d_p)$ where the differential reads
      \beqa
      d_p( u \otimes  x_1 \wedge \ldots \wedge x_p ) = \sum_i (-1)^{i+1}
        u x_i  \otimes x_1 \wedge \ldots \wedge \hat{x}_i\wedge \ldots \wedge x_p   \\
      + \sum_{i<j} (-1)^{i+j} u \otimes   [x_i, x_j]  \wedge  x_1 \wedge \ldots \wedge \hat{x}_i \wedge \ldots
      \wedge \hat{x}_j \wedge \ldots \wedge x_p
      \eeqa
provides a non-minimal projective(in fact free) resolution of $\K$, $C(\g) \stackrel{\epsilon}{\rightarrow} \K $ \cite{W}.
With the latter resolution $C_{\bullet}(\g)$ one 
 calculates homologies of the derived complex $ \K \otimes_{PS}  C_{\bullet}(\g)$ 
 $$E_n=\mathrm{Tor}^{PS}_{n}(\K,\K)  \cong H_{n}(  \K \otimes_{PS}  C_{\bullet}(\g)) = H_{n}(\g, \K)  \ , $$
coinciding with the homologies $H_{n}(\g, \K)$ of the Lie algebra $\g$ with trivial coefficients.
The derived complex $ \K \otimes_{PS}  C_{\bullet}(\g)$ is the chain complex
with degrees ${\bigwedge}^{\bullet} \g=\K\otimes_{PS}PS\otimes {\bigwedge}^{\bullet} \g$
and  differentials $\partial_{p} := id \otimes_{PS} d_p: {\bigwedge}^{p} \g \rightarrow  {\bigwedge}^{p-1} \g $.

The differential $\partial$ is induced by the Lie bracket $[\,\cdot\, ,\,\cdot\,]: {\bigwedge}^2\g \rightarrow \g$ of the graded Lie algebra $\g= \g_1 \oplus \g_2$. It identifies a pair of degree 1 generators $e_i, e_j \in \g_1$ with one degree 2 generator    $e_{ij}:=(e_i \wedge e_j)=[e_i, e_j]\in \g_2$.
The differential $\partial_p$ is a continuation as coderivation(see e.g.\cite{LV}) of the mapping $\partial_2:=-[\,\cdot\, ,\,\cdot\,]$ on the exterior powers ${\bigwedge}^{p} \g$. 
 In greater details the chain degrees read
\beq
   {\bigwedge}^p \g = {\bigwedge}^p \left(V \oplus {\bigwedge}^2 V \right)
   = \bigoplus_{s+r=p} {\bigwedge}^{s} \left({\bigwedge}^2 V \right) \otimes {\bigwedge}^r V
   \eeq
   and   differentials  $\partial_{p=r+s} :
    {\bigwedge}^{s} ({\bigwedge}^2 V) \otimes {\bigwedge}^r(V) \rightarrow {\bigwedge}^{s+1} ({\bigwedge}^2 V) \otimes {\bigwedge}^{r-2}(V) $ 
    are given by
   \beqa
   \label{part}
   \partial_p:\,\,\, e_{i_1 j_1}\wedge \ldots \wedge e_{i_s j_s} \, \, \otimes \,\,  e_1 \wedge \ldots \wedge e_r&\mapsto& \nonumber
   \\
  \sum_{i<j} (-1)^{i+j}  e_{ij} \wedge   e_{i_1 j_1}\wedge \ldots \wedge e_{i_s j_s} &\otimes& e_1 \wedge   \ldots   \wedge \hat{e_i} \wedge \ldots
  \wedge  \hat{e_j} \wedge \ldots \wedge e_r \ .
    \nonumber
   \eeqa
   
By duality, one has the cochain complex $\mathrm{Hom}_{PS}(C(\g),\K)=({\bigwedge}^{\bullet}\g^{\ast},\delta)$
which is a (super)commutative DGA. This cochain complex  calculates the cohomology
 \beq
 \label{iso*}
E_n^{\ast} = \mathrm{Ext}^n_{PS}(\K,\K)\cong H^n(\mathrm{Hom}_{PS}(C(\g),\K)) =H^n(\g,\K) \ .
 \eeq
 The coboundary map $\delta^p:{\bigwedge}^p \g^{\ast}\rightarrow {\bigwedge}^{p+1}\g^{\ast} $ is transposed\footnote{In the presence of metric one has $\delta:=\partial^{\ast}$(see Proposition 1 below).} to the differential $\partial_{p+1}$
 \beqa
 \label{dga}
   \delta^p: \,\,\, e_{i_1 j_1}^{\ast}\wedge \ldots \wedge e_{i_s j_s}^{\ast} \, \, \otimes \,\,  e_1^{\ast} \wedge \ldots \wedge e_r^{\ast}&\mapsto & 
   \\
\sum_{k=1}^{s}  \sum_{i_k<j_k} (-1)^{i+j}    e_{i_1 j_1}^{\ast}\wedge \ldots \wedge \hat{e}_{i_k j_k}^{\ast}\wedge \ldots \wedge e_{i_s j_s}^{\ast} &\otimes&  e_{i_k}^{\ast} \wedge  e_{j_k}^{\ast} \wedge e_1^{\ast} \wedge   \ldots   
  \wedge \ldots \wedge e_r^{\ast} \ ,
   \nonumber
   \eeqa
it is (up to a conventional sign) a continuation  of the dualization of the Lie bracket
$\delta^1:= [\,\cdot\, ,\,\cdot\,]^{\ast}: \g^{\ast}\rightarrow {\bigwedge}^2  \g^{\ast} $ by the Leibniz rule.

 It is important that in the complexes  $({\bigwedge}^{p} \g ,\partial_p)$ and $({\bigwedge}^{p} \g^{\ast} ,\delta^p)$
 two different degees are involved; one is  the homological degree $p:=r+s$ counting the number of
 $\g$-generators, while the second is the tensor degree $t:=2s + r$.
 The differentials $\partial$ and $\delta$  preserve the tensor degree $t$ but the spaces $H_{n}(\g, \K)$ and $H^{n}(\g, \K)$ are not homogeneous in $t$ in general.

 \section{Littlewood formula and $PS$}
 In this section we review the beautiful result of J\'ozefiak and Weyman \cite{JW} giving a
representation-theoretic interpretation of the Littlewood formula
\beq
\label{Littlewood}
\prod_{i} (1- x_i) \prod_{i<j} (1- x_ix_j) = \sum_{\lambda: \lambda=\lambda'} (-1)^{\frac{1}{2}(|\lambda| + r(\lambda))}s_{\lambda} (x)\ .
\eeq
Here the sum is over the self-dual Young diagrams $\lambda$, $s_{\lambda}(x)$ stands for the Schur function and $r(\lambda)$ stands the rank of $\lambda$ which is  the number of diagonal boxes in $\lambda$.

  An irreducible $GL(V)$-module $V_{\lambda}$ is called Schur module,
it has a basis labelled by semistandard Young tableaux which are fillings of the Young diagram $\lambda$
with the numbers of the set $\{1, \ldots ,\dim V \}$.
 The action of the linear group $GL(V)$ on the space $V$ of the generators of the Lie algebra $\g$
 induces a $GL(V)$-action on the UEA $PS=U\g\cong S(V\oplus \Lambda^2 V)$ and  on the 
 space $ {\bigwedge}^{\bullet} \g \cong {\bigwedge}^{\bullet} (V\oplus {\bigwedge}^2 V)$.
  The algebra $PS(V)$ has remarkable property, it is a model of the linear group $GL(V)$, in the sense that it contains every polynomial finite-dimensional 
  irreducible representation $V_{\lambda}$ of $GL(V)$ once and only once
  $$
  PS(V) \cong \bigoplus_{\lambda} V_{\lambda} \ .
  $$
A nice combinatorial proof of this fact was given by  Chaturvedi\cite{Chaturvedi}. The $GL(V)$-model $PS(V)$
enjoys the universal property that every parastatistics Fock representation
specified by the parastatistics order $p\in \mathbb N_0$ is a factor of $PS(V)$
\cite{D-VP},\cite{LP}.


  The differential $\partial$ 
   commutes 
  with the $GL(V)$ action and  the homology  $ H_{\bullet}(\g,\K)$ is also a $GL(V)$-module.
The decomposition of the $GL(V)$-module $H_n(\g, \K)$ into irreducible polynomial representations $V_{\lambda}$ is given by the following theorem;
 \begin{theorem}[J\'ozefiak and Weyman\cite{JW}, Sigg\cite{Sigg}]
 \label{jwthm} The homology  $ H_{\bullet}(\g,\K)$ of the 2-nilpotent Lie algebra $\g=V \oplus {\bigwedge}^2 V$
 decomposes into irreducible  $GL(V)$-modules 
  \beq
 H_n(\g, \K)= H_n({\bigwedge}^{\bullet} \g, \partial ) \cong \mathrm{Tor}^{PS}_n(\K,\K)(V)\cong\bigoplus_{\lambda:\lambda = \lambda' } V_{\lambda}
  \eeq
 where  the sum is over self-conjugate Young diagrams $\lambda$ such that $n = \frac{1}{2}(|\lambda| + r(\lambda))$. 
  \end{theorem}
  The data 
  $H_n(\g, \K)=\mathrm{Tor}^{PS}_n(\K,\K)$ encodes 
the minimal free resolution $P_{\bullet}$ (cf. (\ref{minP})). 

 The acyclicity of the complex $P_{\bullet}$ implies an identity
about  the $GL(V)$-characters 
$$
ch \,  PS(V) \, . \, ch \left(\bigoplus_{\lambda:\lambda = \lambda' } (-1)^{\frac{1}{2}(|\lambda| + r(\lambda))} V_{\lambda}\right) =1 \ .
$$
The character of a Schur module $V_{\lambda}$ is the Schur function,
  $ch  \, V_{\lambda}= s_{\lambda}(x)$. Due to  the Poincar\'e-Birkhoff-Witt theorem
  $ch \,  PS(V)= ch \, S(V \oplus {\bigwedge}^2 V)$ thus the identity reads
\beq
 \prod_i \frac{1}{(1- x_i)} \prod_{i<j}\frac{1}{(1- x_i x_j) }  \sum_{\lambda: \lambda=\lambda'} (-1)^{\frac{1}{2}(|\lambda| + r(\lambda) ) } s_{\lambda} (x)=1 \ .
  \eeq
But the latter identity is nothing but rewriting of the Littlewood identity (\ref{Littlewood}).
The moral is that the Littlewood identity reflects a homological property of the algebra $PS$, namely
the above particular structure  of the minimal projective (free) resolution of $\K$ by $PS$-modules. 
 
\section{Homotopy algebras $A_{\infty}$ and $C_{\infty}$}
\begin{definition}{($A_{\infty}$-algebra)} A homotopy associative  algebra, 
or  $A_{\infty}$-algebra,
over $\K$ is a $\Z$-graded vector space
$
  A = \bigoplus_{i\in \Z} A^i \nonumber
  $
endowed with a family of graded 
mappings (operations)
\beq
  m_n : A^{\otimes n} \rightarrow A ,  \qquad  \deg  (m_n)=2-n \qquad n \geq 1 \nonumber
\eeq
satisfying  the Stasheff identities  $\bf SI(n)$  for   $n\geq 1$
\beq
  \sum_{r+s+t=n} (-1)^{r+st} m_{r+1+t}(Id^{\otimes r}\otimes m_s \otimes Id^{\otimes t}) = 0
  \qquad \bf SI(n) \nonumber
  \eeq
 where the   sum runs over all decompositions  $n=r+s+t$.
\end{definition}
Here we assume the Koszul sign convention $(f\otimes g)(x \otimes y)= (-1)^{|g||x|} f(x)\otimes g(y)$.
We define the shuffle product $Sh_{p,q}: A^{\otimes p}\otimes A^{\otimes q} \rightarrow A^{\otimes p+q}$ throughout the expression
$$
(a_1\otimes \ldots \otimes a_p) \shuffle (a_{p+1}\otimes \ldots \otimes a_{p+q})=
\sum_{\sigma\in Sh_{p,q}} {sgn(\sigma)}\, a_{\sigma^{-1}(1)}\otimes \ldots \otimes a_{\sigma^{-1}(p+q)}
$$
where the sum runs over all $(p,q)$-shuffles $Sh_{p,q}$, i.e., over all permutations $\sigma\in S_{p+q}$
such that 
$\mbox{\small $\sigma(1) < \sigma(2)< \ldots < \sigma(p)$}$ and $
\mbox{\small$\sigma(p+1) < \sigma(p+2)< \ldots < \sigma(p+q)$} \ . $
\begin{definition}{($C_{\infty}$-algebra \cite{Kadeishvili})} A homotopy commutative  algebra,
or $C_{\infty}$-algebra,
 is an $A_{\infty}$-algebra $\{A, m_n \}$
with the 
 condition: each operation $m_n$ vanishes
on shuffles
\beq
m_n \left((a_1\otimes \ldots \otimes a_p) \shuffle (a_{p+1}\otimes \ldots \otimes a_n)\right) = 0\ ,
\qquad 1 \leq p \leq n-1 \ .
\eeq
\end{definition}  
In particular for $m_2$ we have $m_2(a\otimes b \pm b\otimes a)=0$, so a $C_\infty$-algebra
 such that $m_n=0$ for $n\geq3$ is a (super)commutative DGA.
 
A morphism of two $A_{\infty}$-algebras  $A$ and $B$ is a family of graded maps 
$ f_n: A^{\otimes n} \rightarrow B $ for  $n\geq 1$  with $\deg{f_n}=1-n$  
such that 
the following conditions hold
$$\sum_{r+s+t=n} (-1)^{r+st} f_{r+1+t} (Id^{\otimes r} \otimes m_s \otimes Id^{\otimes r}) =
  \sum_{1\leq q \leq n} (-1)^{S} m_q(f_{i_1} \otimes f_{i_2} \otimes \ldots \otimes f_{i_q})
$$
where the sum on the RHS is over all decompositions $i_1+ \ldots + i_q=n$ and the sign  is determined
by $S= \sum_{k=1}^{q-1}(r-q)(i_q - 1) $. 
The morphism $f$ is a {\it quasi-isomorphism of $A_{\infty}$-algebras} if $f_1$ is a quasi-isomorphism.
It is strict if $f_i=0$ for all $i\neq 1$. The identity morphism of $A$
 is the strict morphism $f$ such that
$f_1$ is the identity of $A$.

A morphism of $C_{\infty}$-algebras is a morphism of $A_{\infty}$-algebras with components
vanishing on shuffles
$
f_n \left((a_1\otimes \ldots \otimes a_p) \shuffle (a_{p+1}\otimes \ldots \otimes a_n)\right) = 0\ ,
\,\, 1 \leq p \leq n-1 \ .
$

\section{Homotopy Transfer Theorem}
\begin{lemma}[see e.g.\cite{LV}] Every cochain complex $(A,d)$ of vector spaces over a field $\K$
has its cohomology $H^{\bullet}(A)$ as a deformation retract.
\end{lemma}
One can always choose a vector space decomposition of the 
cochain complex $(A,d)$ such that $A^n\cong B^n \oplus H^n \oplus B^{n+1}$ where $H^n$ is the cohomology and $B^n$ is the space of coboundaries, $B^{n}=d A^{n-1}$. We choose a homotopy $h:A^n \rightarrow A^{n-1}$ which  identifies $B^{n}$ with its copy in $A^{n-1}$ and is 0 on $H^n \oplus B^{n+1}$.
The projection $p$ to the cohomology  and the cocycle-choosing inclusion $i$ given by
$\xymatrix{ A^n    \ar@<.5ex>[r]^{p} &  H^{n}\ar@<.5ex>[l]^{i}} $ are chain homomorphisms
(satisfying the additional conditions $hh=0$, $hi=0$ and $ph=0$).
With these choices done the complex $(H^{\bullet}(A),0)$ is a
deformation retract of $(A,d)$ 
\[
\xymatrix{
 \ar@(ld,lu)^{h}}\!\!\!\!\!
\xymatrix{
(A, d)    \ar@<.5ex>[r]^{p \hspace{10pt}   }
&  (H^{\bullet}(A),0) \ar@<.5ex>[l]^{i  \hspace{10pt}  }
 } \ , \qquad
pi = Id_{H^{\bullet}(A)} \ , \qquad ip - Id_A = d h + h d  \ .\]

Let now $(A,d,\mu)$ be a  DGA, i.e.,  $A$ is endowed with an associative product $\mu$ compatible with $d$.
 The cochain complexes $(A,d)$ and its contraction  $H^{\bullet}(A)$ 
 are homotopy equivalent, but the associative 
 structure
 is not
 stable under homotopy equivalence. However the associative structure on
 $A$ can be transferred to an $A_{\infty}$-structure on a homotopy equivalent complex,
 a particular interesting complex being the deformation retract $H^{\bullet}(A)$. 
 For a friendly introduction to homotopy transfer theorems 
in much boarder context we send the reader to the textbook \cite{LV}.

  \begin{theorem}[Kadeishvili\cite{Kadeishvili}]
  \label{cinf}
  Let $(A,  d, \mu)$ be a (commutative) DGA over a field $\K$.
  There exists a $A_{\infty}$-algebra ($C_{\infty}$-algebra) structure on the cohomology  $H^{\bullet}(A)$ and a $A_{\infty}(C_{\infty})$-quasi-isomorphism
  $f_i:({\otimes}^i H^{\bullet}(A), \left\{ m_i \right\})\rightarrow (A,\left\{d, \mu, 0, 0, \ldots \right\} )$
   such that the inclusion $f_1=i:H^{\bullet}(A)\rightarrow A$ is a cocycle-choosing homomorphism
  of  cochain complexes. The differential on $H^{\bullet}(A)$  is zero $m_1=0$ and $m_2$ is strictly associative operation  induced by the multiplication  on $A$. The resulting structure is unique up to quasi-isomorphism. \end{theorem}

   Kontsevich and Soibelman\cite{KS} gave an explicit expressions for the higher operations of the induced
  $A_\infty$-structure  as  sums over decorated  planar binary  trees with one root where all leaves are decorated by the 
  inclusion $i$, the root by the projection $p$ the vertices by the  product $\mu$ of the (commutative) DGA $(A,d,\mu)$ and the internal edges by the homotopy $h$. The $C_\infty$-structure implies additional symmetries on trees.
  
 \noindent
We will make use of the  graphic representation for the  binary operation on $H^{\bullet}(A)$
   $$\xymatrix{& \ar[dr]_i& &  \ar[dl]^i \\
         m_2(x,y):= p \mu(i(x),i(y)) \quad \mbox{or} \quad     m_2=&&\mu \ar[d]_p& \\
              &&&} $$
              and the ternary one $m_3(x,y,z) = p \mu(i(x), h \mu(i(y),i(z)))- p \mu( h \mu(i(x),i(y)),i(z))$ being the sum of  two planar binary trees with three leaves 
              $$  \ba{ccc}
\xymatrix{ \ar[ddrr]_i& &\ar[dr]_i&&\ar[dl]^i\\
              &&&\mu \ar[dl]^h &\\
            m_3= \quad &&\mu \ar[d]_p &&  \\
              && &&} & \ba{c} \vspace{-5cm}  - \\ \ea&
             \xymatrix{\ar[dr]_i & &\ar[dl]^i&&\ar[ddll]_i\\
              &\mu \ar[dr]^h && &\\
              &&\mu \ar[d]_p && \\
              && &&}\ea$$
  \begin{proposition} The cohomology
 $H^{\bullet}(\g,\K)\cong \mathrm{Ext}^{\bullet}_{PS}(\K,\K)$ of the 2-nilpotent graded Lie algebra $\g=V\otimes {\bigwedge}^2 V$
  is a homotopy commutative algebra. The
$C_{\infty}$-algebra $H^{\bullet}(\g,\K)$ is  generated 
in degree 1, i.e., in $H^{1}(\g,\K)$  by the operations $m_2$ and $m_3$.
 \end{proposition}
 {\bf Sketch of the proof.} 
We apply  the Kadeishvili  homotopy transfer Theorem \ref{cinf} for the  commutative DGA $({\bigwedge}^{\bullet} \g^{\ast},\mu, \delta^{\bullet})$ and its deformation retract $H^{\bullet}({\bigwedge}^{\bullet} \g^{\ast})\cong H^{\bullet}(\g,\K)$ and conclude that the cohomology  $H^{\bullet}(\g,\K)$ is a $C_{\infty}$-algebra.

 Further on we need
  the  explicit mappings in the deformation retract.
 Let us choose a metric $g(\,\cdot\, ,\,\cdot\,)=\langle \,\cdot\, , \,\cdot\,\rangle $
   on the vector space $V$
  and an orthonormal basis 
   $\langle  e_i ,  e_j\rangle= \delta_{ij}$. The choice  induces  a metric on 
    ${\bigwedge}^{\bullet} \g \stackrel{g}{\cong} {\bigwedge}^{\bullet} \g^{\ast}$.
    In the presence of metric $g$
   the differential $\delta$ 
 is   identified with the adjoint of  $\partial$, $\delta\stackrel{g}{:=}\partial^{\ast}$
   while $\partial$ plays the role of a homotopy. Then the deformation  retract $H^{\bullet}({\bigwedge}^{\bullet} \g^{\ast}, \delta^{\bullet})$  of $({\bigwedge}^{\bullet} \g^{\ast}, \delta^{\bullet})$
    looks like
    $$
pi=Id_{H^{\bullet}({\bigwedge}^{\bullet} \g^{\ast})} \ ,\quad \qquad ip - Id_{{\bigwedge}^{\bullet} \g^{\ast}} = \delta \delta^{\ast} + \delta^{\ast} \delta \ ,\qquad \delta^{\ast}\stackrel{g}{=}\partial \ .
   $$
Here the projection $p$ identifies the subspace $\ker \delta \cap \ker \delta^{\ast}$ with $H^{\bullet}({\bigwedge}^{\bullet} \g^{\ast})$,
which is the orthogonal complement of the space of the coboundaries ${\rm im} \delta$.
The cocycle-choosing homomorphism $i$  is $Id$ on  $H^{\bullet}({\bigwedge}^{\bullet} \g^{\ast})$ and zero on coboundaries. 

 Due to  the isomorphisms $\mathrm{Tor}^{PS}_{n}(\K,\K) \cong {\rm Ext}^n_{PS}(\K,\K)$ (see eq. (\ref{cartaniso}))  
 induced by  $V\stackrel{g}{\cong} V^{\ast}$  the theorem \ref{jwthm} implies the decomposition of
    $H^{\bullet}(\g,\K)$ into  Schur modules 
    \beq
   H^{n}(\g,\K)\cong H^{n}({\bigwedge} \g^{\ast}, \delta)\cong \mathrm{Ext}_{PS}^{n}(\K,\K)(V^{\ast}) \cong{\bigoplus}_{\lambda:\lambda = \lambda' } V_{\lambda}
\nonumber  \eeq
 where  the sum is over self-conjugate diagrams  $\lambda$ such that $n = \frac{1}{2}(|\lambda| + r(\lambda))$. 
The operation $m_n$ is bigraded  $\deg_{n',t'}(m_n)=(2-n, 0)$  by homological degree $n'$ and tensor degree $t'$(weight).
 The bi-grading impose the vanishing of many higher  products.

The Kontsevich and Soibelman tree representations of the operations $m_n$ provide explicit 
expressions.
Let us take $\mu$ to be the super-commutative product $\wedge$ on the DGA $({\bigwedge}^{\bullet} \g^{\ast}, \delta^{\bullet})$.
The projection  $p$ maps  onto the Schur modules $V_{\lambda}$ with $\lambda =\lambda'$.
   
The binary operation on the degree 1 generators $e_i\in H^1(\g,\K)$
is trivial, one gets $$m_2(e_i,e_j)= p  (e_i\wedge e_j)=0 \qquad  p(V_{(1^2)})=0.$$
Hence $H^{\bullet}(\g,\K)$ could not be generated in $H^{1}(\g,\K)$ as  algebra with product $m_2$.
 
 The ternary operation $m_3$ restricted to $H^1(\g,\K)$ is nontrivial, indeed one has
 \beqa
 m_3(e_i, e_j ,e_k)&= &p\left\{ e_i\wedge \partial(e_j\wedge e_k) - \partial(e_i\wedge e_j)\wedge e_k\right\}=
  p\left\{e_{ij}\wedge e_k -e_i\wedge e_{jk}  \right\} \nonumber\\
 &=& p\left\{ (e_{ij}\wedge e_{k} +e_{jk} \wedge e_{i} +e_{ki} \wedge e_{j})-e_{ki} \wedge e_{j}\right\}=
  e_{ik} \wedge e_{j} \in H^{2}(\g,\K) \nonumber
 \eeqa
The completely antisymmetric combination in the brackets $(\ldots)$ spans the Schur module $V_{(1^3)}$, 
$p(e_{ij}\wedge e_{k} +e_{jk} \wedge e_{i} +e_{ki} \wedge e_{j})=0$ yields a Jacobi-type identity.  The monomials $e_{ij} \wedge e_{k}$ modulo $V_{(1^3)}$ span a Schur module $V_{(2,1)}\in H^{2}(\g,\K)$ with basis  in bijection with the semistandard  Young tableaux 
$
 e_{ik} \wedge e_j \leftrightarrow { \footnotesize \ba{|c|c|}\hline
\,\, i \, \,\, & \, j \,\, \\ \hline
 k \\
\cline{1-1}
\ea} $ and $e_{ij}\wedge e_k  \leftrightarrow {\footnotesize \ba{|c|c|}\hline
\,\, i \,\,\,&\, k\,\,  \\ \hline
 j \\
\cline{1-1}
\ea}$.

 We check  the symmetry condition on ternary operation $m_3$ in $C_{\infty}$-algebra;
 indeed $m_3$ vanishes on the (signed) shuffles $Sh_{1,2}$ and $Sh_{2,1}$
 $$ m_3(e_i\shuffle e_j\otimes e_k)=m_3(e_i, e_j, e_k) -m_3(e_j, e_i, e_k)+m_3( e_j, e_k,e_i) = 0=m_3(e_i\otimes e_j\shuffle e_k). $$

On the level of Schur modules the ternary operation glues  three fundamental $GL(V)$-representations $V_\Box$ into a Schur module $V_{(2,1)}$. By  iteration of the process of gluing boxes we generate all  elementary hooks $V_k:=V_{(1^k, k+1)}$, 
\[
m_3(V_\Box,V_\Box,V_\Box)=  V_{\mbox{\tiny
$\ba{|c|c|}\hline
\hspace{3pt} &  \hspace{3pt} \\ \hline
 \\
\cline{1-1}
\ea$}} \, \,, \quad  m_3\left( V_\Box, V_{\mbox{\tiny
$\ba{|c|c|}\hline
\hspace{3pt} &  \hspace{3pt} \\ \hline
 \\
\cline{1-1}
\ea$}} \, , V_\Box\right)= 
V_{\mbox{\tiny$\ba{|c|c|c|}\hline
\hspace{3pt} &  \hspace{3pt} & \hspace{3pt} \\ \hline
 \\
\cline{1-1}\\
\cline{1-1}
\ea$}} \, \, ,  \ldots ,  
m_3(V_0, V_{k}, V_0)=V_{k+1} \ .
\]
In our context the more convenient notation for Young diagrams is due to Frobenius: $\lambda:=(a_1, \ldots,a_r|b_1, \ldots b_r)$ stands 
for a diagram $\lambda$ with $a_i$ boxes  in the $i$-th row on the right of the diagonal,
and with $b_i$  boxes  in the $i$-th coloumn below the diagonal  and the rank $r=r(\lambda)$
is the number of boxes on the diagonal. 

For self-dual diagrams $\lambda=\lambda'$, i.e., $a_i=b_i$
 we set
  $V_{a_1, \ldots, a_r}:=V_{(a_1, \ldots,a_r|a_1, \ldots a_r)}$
when  $a_1>a_2> \ldots >a_r \geq 0$ (and set the convention $V_{a_1, \ldots, a_r}:=0$ otherwise).
Any two elementary hooks $V_{a_1}$ and $V_{a_2}$  can be glued together by the binary  operation $m_2$,
the decomposition of $m_2(V_{a_1},V_{a_2})\cong m_2(V_{a_2},V_{a_1})$ is given by
\[
 m_2(V_{a_1},V_{a_2})=V_{{a_1},{a_2}}\oplus (\bigoplus_{i=1}^{{a_2}} V_{{a_1}+i, {a_2}-i}) \qquad {a_1}\geq {a_2}
\]
where the ``leading'' term $V_{{a_1},{a_2}}$ has the diagram with minimal height. 
Hence any $m_2$-bracketing of the hooks $V_{a_1},V_{a_2} \ldots , V_{a_r}$ yields\footnote{The operation $m_2$ 
is  associative 
thus the result does not depend on the 
choice of the bracketing.} 
a sum of $GL(V)$-modules 
$$
 m_2(\ldots m_2(m_2(V_{a_1},V_{a_2}),V_{a_3}), \ldots , V_{a_r})=
V_{a_1, \ldots, a_r} \oplus \ldots
$$
whose module with  minimal height   is precisely   $V_{a_1, \ldots, a_r}$. We conclude that all 
elements in the $C_{\infty}$-algebra $H^{\bullet}(\g,\K)$ can be generated in $H^{1}(\g,\K)$ by
$m_2$ and $m_3$. $\Box$

\begin{acknowledgement}
We are grateful to Jean-Louis Loday for many enlightening discussions and his 
encouraging interest. The work was  supported by the French-Bulgarian Project Rila under the
contract Egide-Rila N112.
\end{acknowledgement}
\end{document}